\documentclass[preprint,review,10pt]{elsarticle}
\usepackage{amsfonts}
\usepackage{amsmath,amssymb}
\usepackage{graphicx}
\usepackage{setspace}
\usepackage[left=2cm,top=2cm,right=2cm]{geometry}

\newtheorem{definition}{Definition}

\newtheorem{lemma}{Lemma}

\newtheorem{theorem}{Theorem}
\numberwithin{equation}{section} \journal{XYZ}
\begin{document}
\title{On Generalized Sz\'{a}sz-Mirakyan operators}
\author[label*,label1,label2]{Prashantkumar Patel}
\ead{prashant225@gmail.com}
\author[label1,label3]{Vishnu Narayan Mishra}
\ead{vishnu\_narayanmishra@yahoo.co.in; vishnunarayanmishra@gmail.com}
\address[label1]{Department of Applied Mathematics \& Humanities,
S. V. National Institute of Technology, Surat-395 007 (Gujarat), India}
\address[label2]{Department of Mathematics,
St. Xavier's College(Autonomous), Ahmedabad-380 009 (Gujarat), India}
\address[label3]{L. 1627 Awadh Puri Colony Beniganj, Phase-III, Opposite - Industrial Training Institute (ITI), Ayodhya Main Road, Faizabad, Uttar Pradesh
224 001, India}
\fntext[label*]{Corresponding authors}
\begin{abstract}
In the present article, we propose the generalization of  Sz\'{a}sz-Mirakyan operators, which is a class of linear positive operators of discrete type depending on a real parameters. We give theorem of degree of approximation, the Voronovskaya Asymptotic formula and statistical convergence.
\end{abstract}

\maketitle
\section{Introduction}
Since 1912, Bernstein Polynomial and its various generalization have been studied by Bernstein \cite{Berstien1912Demo}, Sz\'{a}sz \cite{szasz1950generalization}, Meyer-Konig and Zeller \cite{meyer1960bernsteinsche}, Cheney and Sharma \cite{cheney1964generalization}, Stancu \cite{stancu1968approximation}. Bernstein polynomials are based on binomial and negative binomial distributions. In 1941, Sz\'{a}sz and Mirakyan \cite{mirakyan1941approximation} have introduced operator using the Poisson distribution.  After that various generalization and approximation properties was discussed by many researcher of the Sz\'{a}sz-Mirakyan operators. We mention that rate of convergence develop by Rempulska and Walczak \cite{rempulska2001approximation}, asymptotic expansion introduced by Abel \textit{et al.} \cite{abel2007asymptotic}. In 1976, May \cite{may1976saturation} showed that the Baskakov operators can reduce to the Sz\'{a}sz-Mirakyan operators. In 1972, Jain introduced another generalization of Sz\'{a}sz-Mirakyan operators in \cite{jain1972approximation}  and discussed the relation between the local smoothness of function and local approximation, the degree of approximation and the statistical convergence of the Jain operators was studied by Agratini \cite{agratini2013approximation}. Umar and Razi \cite{umar1985approximation} studied Kantorovich-type extension of Jain operators.
Durrmeyer type generalization of Jain operators and its approximation properties was elaborated by   Tarabie \cite{tarabie2012jain}, Mishra and Patel \cite{mishrasome2013}, Patel and Mishra \cite{Patelmishra20131} and Agratini \cite{agratini2014approximation}. See some related work in this area \cite{Deepmala2014,mishra2013inverse,mishra2013statistical}.  Motivated by this we give further modification of Jain operators in this paper with the help of a generalized Poisson type distribution, consider its convergence properties and give its degree of approximation, asymptotic formula, the statistical convergence. The results for the Sz\'{a}sz-Mirakyan operator and Jain operators can easily be obtained from our operator as a particular case.\\
In the year 1902, Jensen \cite{jensen1902identite} use the following Lagrange's formula
\begin{equation}\label{14.1.eq1}
\phi(z) = \phi(0) + \sum_{k=1}^{\infty} \frac{1}{k!}\left[\frac{d^{k-1}}{dz^{k-1}}((f(z)^k)\phi'(z)\right]_{z=0} \left(\frac{z}{f(z)}\right)^k.
\end{equation}
Using equation (\ref{14.1.eq1}) Jain \cite{jain1972approximation} discussed convergence of Jain operators. Here we use above Lagrange's formula of discuss convergence of generalized Jain operators. Proceed by setting
$$ \phi(z) = a^{\alpha z}~~~~~  and  ~~~~~f(z) = a^{\beta z}, ~~~~1<a\leq e; ~ 0<\alpha<\infty;~-1<\beta<1;$$
we shall get
\begin{eqnarray}\label{14.1.eq2}
a^{\alpha z} &=& 1 + \alpha\sum_{k=1}^{\infty} \frac{1}{k!} (\log a)^k (\alpha + \beta k )^{k-1}  \left(\frac{z}{a^{\beta z}}\right)^k\nonumber\\
&=& \sum_{k=0 }^{\infty} \alpha (\log a)^k (\alpha + \beta k )^{k-1}  \frac{u^k}{k!},~~ u= za^{-\beta z},
\end{eqnarray}
where $z$ and $u$ are sufficiently small such that $|\beta u | < a^{-1}$ and $| \beta z|<1$.
By taking z = 1, we have
\begin{equation}\label{14.1.eq3}
1  = \sum_{k=0 }^{\infty} \frac{\alpha}{k!} (\log a)^k (\alpha + \beta k )^{k-1}a^{-(\alpha + \beta k)}.
\end{equation}
From above discussion proof the following Lemma is achieve.
\begin{lemma}
For $1<a\leq e$; $0<\alpha<\infty$ and $|\beta| <1$, let
\begin{equation}\label{14.1.eq4}
\omega_{\beta,a}(k,\alpha) = \alpha(\log a)^k ( \alpha + k\beta)^{k-1} \frac{a^{-(\alpha+k\beta)}}{k!}
\end{equation} then\vspace{-0.5cm}
\begin{equation}\label{14.1.eq5}
\sum_{k=0}^{\infty} \omega_{\beta,a}(k,\alpha)= 1.
\end{equation}
\end{lemma}
\begin{lemma}
Let
\begin{equation}\label{14.1.eq6}
 S(r,\alpha,\beta,a) = \sum_{k=0 }^{\infty} \frac{1}{k!} (\log a)^k (\alpha + \beta k )^{k+r-1}a^{-(\alpha + \beta k)}
 \end{equation}
and\vspace{-0.5cm}
\begin{equation}\label{14.1.eq7}\alpha S(0,\alpha,\beta,a) = 1. \end{equation}
Then
\begin{equation}\label{14.1.eq8}
S(r,\alpha,\beta,a) = \sum_{k=0}^{\infty} (\beta \log a)^k (\alpha + k\beta) S( r-1, \alpha+k\beta,\beta).
\end{equation}
\end{lemma}
\textbf{Proof:}
It can easily be seen that,
\begin{eqnarray}\label{14.1.eq9}
 S(r,\alpha,\beta,a) &=& \sum_{k=0 }^{\infty} \frac{1}{k!} (\log a)^k (\alpha + \beta k )(\alpha + \beta k )^{k+r-2}a^{-(\alpha + \beta k)}\nonumber\\
 &=& \alpha S(r-1,\alpha,\beta) + \beta\sum_{k=1 }^{\infty} \frac{\alpha}{(k-1)!} (\log a)^k    (\alpha + \beta k )^{k+r-2}a^{-(\alpha + \beta k)}\nonumber\\
 &=& \alpha S(r-1,\alpha,\beta) + \beta \log a S(r,\alpha+\beta,\beta).
\end{eqnarray}
Repeated use of equation (\ref{14.1.eq9}), we get
\begin{eqnarray*}
S(r,\alpha,\beta,a) &=& \alpha S(r-1,\alpha,\beta,a ) + \beta \log a \left[(\alpha+\beta) S(r-1,\alpha+\beta ,\beta,a) + \beta \log a S(r,\alpha+2\beta,\beta,a)\right]\\
&=& \alpha S(r-1,\alpha,\beta,a) + \beta \log a (\alpha+\beta) S(r-1,\alpha+\beta ,\beta,a) + (\beta \log a)^2 S(r,\alpha+2\beta,\beta,a)\\
&=& \sum_{k=0}^{\infty} (\beta \log a)^k (\alpha + k\beta)S(r-1,\alpha+k\beta,\beta,a),
\end{eqnarray*}
which is required result.\\
Now, when $|\beta \log a | < 1$ , we have
\begin{equation}\label{14.1.eq10}
S(1,\alpha,\beta,a) =\sum_{k=0}^{\infty} (\beta \log a)^k = \frac{1}{1- \beta \log a};
\end{equation}
\begin{equation}\label{14.1.eq11}
S(2,\alpha,\beta,a) = \sum_{k=0}^{\infty} \frac{(\beta \log a)^k(\alpha+k \beta)}{1-\beta \log a} = \frac{\alpha}{(1-\beta \log a)^2}+ \frac{\beta^2 \log a }{(1-\beta \log a)^3};
\end{equation}
\begin{eqnarray}\label{14.1.eq12}
S(3,\alpha,\beta,a) &=& \sum_{k=0}^{\infty}(\beta \log a)^k(\alpha+k \beta) S(2,\alpha+k\beta,\beta,a)\nonumber\\
 &=& \frac{1}{(1-\beta \log a)^2}\sum_{k=0}^{\infty}(\beta \log a)^k(\alpha+k \beta)^2+ \frac{\beta^2 \log a }{(1-\beta \log a)^3}\sum_{k=0}^{\infty}(\beta \log a)^k(\alpha+k \beta) \nonumber\\
&=&\frac{1}{(1-\beta \log a)^2}\left[ \frac{\alpha^2}{(1-\beta \log a)} + \frac{2 \alpha\beta^2 \log a }{(1-\beta \log a)^2} + \frac{\beta^3 \log a(1+\beta\log a ) }{(1-\beta \log a)^3}\right]\nonumber\\
&& +  \frac{ \beta^2 \log a }{(1-\beta \log a)^3} \left[\frac{\alpha}{1-\beta\log a } + \frac{\beta^2 \log a }{(1-\beta\log a)^2}\right]\nonumber\\
 &=& \frac{\alpha^2}{(1-\beta \log a)^3} + \frac{3 \alpha\beta^2\log a }{(1-\beta \log a)^4} + \frac{ (\beta^3 + 2\beta^4)\log a  }{(1-\beta \log a)^5}
\end{eqnarray}
and
\begin{eqnarray}\label{14.1.eq13}
S(4,\alpha,\beta,a) &=& \sum_{k=0}^{\infty}(\beta \log a)^k(\alpha+k \beta) S(3,\alpha+k\beta,\beta,a)\nonumber\\
 &=& \sum_{k=0}^{\infty}(\beta \log a)^k(\alpha+k \beta)\left[\frac{(\alpha+k\beta)^2}{(1-\beta \log a)^3} + \frac{3 (\alpha+k\beta)\beta^2\log a }{(1-\beta \log a)^4} + \frac{ (\beta^3 + 2\beta^4)\log a  }{(1-\beta \log a)^5}\right] \nonumber\\
 &=& \frac{1}{(1-\beta \log a)^3} \sum_{k=0}^{\infty}(\beta \log a)^k(\alpha+k \beta)^3 + \frac{3\beta^2\log a }{(1-\beta \log a)^4}\sum_{k=0}^{\infty}(\beta \log a)^k(\alpha+k \beta)^2 \nonumber\\
&& + \frac{ (\beta^3 + 2\beta^4)\log a  }{(1-\beta \log a)^5} \sum_{k=0}^{\infty}(\beta \log a)^k(\alpha+k \beta)\nonumber\\
&=& \frac{\alpha^3 }{(1-\beta \log a)^4} +  \frac{6 \alpha^2 \beta^2\log a }{(1-\beta \log a)^5 } +  \frac{2\alpha\beta^3(2 +\beta)\log a+ 9\alpha \beta^4(\log a)^2 }{(1-\beta \log a)^6 }\nonumber\\
&&+ \frac{\beta^4 \log a +2\beta^5 (4  +\beta) (\log a)^2 + 4\beta^6(\log a)^3}{(1-\beta \log a)^7 }.
\end{eqnarray}
\begin{definition}
We may now define the operator as
\begin{equation}\label{14.1.eq14}
 P_{n}^{[\beta,~a]}(f,x) =  \sum_{k=0}^{\infty} \omega_{\beta,a}(k,n x) f\left(\frac{k}{n}\right)
 \end{equation}
where $f\in C(\mathbf{R}^{+})$; $0\leq \beta < 1$; $1< a \leq e$ and $\omega_{\beta,a}(k,n x)$ as defined in equation (\ref{14.1.eq4}).
\end{definition}
 In the particular case, $a=e$ the operators defined in equation (\ref{14.1.eq14}) convert to Jain operator \cite{jain1972approximation};  $\beta=0$, $a=e$, the operator $P_{n}^{[\beta,~a]}$, $n\in \mathbf{N}$, turns into well-known Sz\`{a}sz-Mirakjan operators \cite{szasz1950generalization}. This operators has different approximation properties then the operators discussed in \cite{szasz1950generalization,jain1972approximation}. \vspace{-0.3cm}
 \section{Estimation of Moments}\vspace{-0.5cm}
 We required following results to prove main results.\vspace{-0.4cm}
 \begin{lemma}The operators $P_{n}^{[\beta,~a]}$ $n \geq 1$, defined by (\ref{14.1.eq14}) satisfy the following relations
 \begin{enumerate}
 \item $P_{n}^{[\beta,a]}(1,x)=1;$
 \item $\displaystyle P_{n}^{[\beta,~a]}(t,x) = \frac{x\log a}{1- \beta \log a};$
 \item $\displaystyle P_{n}^{[\beta,~a]}(t^2,x) = \frac{x^2(\log a)^2}{(1-\beta \log a)^2} +\frac{x \log a}{n(1-\beta  \log a)^3};$
 \item $\displaystyle P_{n}^{[\beta,~a]}(t^3,x) = \frac{x^3(\log a)^3 }{(1-\beta \log a)^3} + \frac{ 3 x^2 (\log a)^2 }{n (1-\beta \log a)^4} + \frac{x\log a(1+2 \beta  \log a +2 \beta ^4 (\log a)^3 -2 \beta ^4 (\log a)^4)}{n^2(1-\beta \log a)^5};$
 \item $\displaystyle P_{n}^{[\beta,~a]}(t^4,x) = \frac{x^4 (\log a) ^4}{(1-\beta  \log a )^4} + \frac{6 x^3 (\log a) ^3}{n (1-\beta  \log a )^5} + \frac{x^2(\log a)^2 \left(7 +8 \beta  \log a+2 \beta ^4 (\log a) ^3-2 \beta ^4 (\log a) ^4\right)}{n^2 (1-\beta  \log a )^6}\\
     ~~~~~~~~~~~~~~~~~+\frac{x \left((\log a) +8 \beta  (\log a) ^2+6 \beta ^2 (\log a) ^3+(12 \beta ^4 (\log a)^4-16 \beta ^5 (\log a) ^5+6 \beta ^6 (\log a) ^6)(1-\log a) \right)}{n^3 (1-\beta  \log a )^7}.$
      \end{enumerate}
 \end{lemma}
 \textbf{Proof:}
 By the relation (\ref{14.1.eq3}), it clear that $P_{n}^{[\beta,a]}(1,x)=1$.\\
 By the simple computation, we get
 \begin{eqnarray*}
P_{n}^{[\beta,~a]}(t,x) &=& n x \sum_{k=0}^{\infty} (\log a)^k ( n x + k\beta)^{k-1} \frac{a^{-(nx+k\beta)}}{k!} \frac{k}{n}
= x\log a S(1,n x+\beta,\beta,a)= \frac{x\log a}{1- \beta \log a}.\\
P_{n}^{[\beta,~a]}(t^2,x) &=& n x \sum_{k=0}^{\infty} (\log a)^k ( n x + k\beta)^{k-1} \frac{a^{-(nx+k\beta)}}{k!} \frac{k^2}{n^2}\\
&=& \frac{x}{n}\left\{ (\log a)^2S(2,n x+2\beta,\beta,a) + \log a S(1,n x+\beta,\beta,a)\right\}\\
&=& \frac{x^2(\log a)^2}{(1-\beta \log a)^2} +\frac{x \log a}{n(1-\beta  \log a)^3}.\\
P_{n}^{[\beta,~a]}(t^3,x) &=& n x \sum_{k=0}^{\infty} (\log a)^k ( n x + k\beta)^{k-1} \frac{a^{-(nx+k\beta)}}{k!} \frac{k^3}{n^3}\\
&=& \frac{x\log a }{n^2}\left\{(\log a)^2 S(3,n x+3\beta,\beta,a) + 3\log a S(2,n x+2\beta,\beta,a) +S(1,n x+\beta,\beta,a) \right\}\\
&=& \frac{x\log a }{n^2}\left\{(\log a)^2 \left( \frac{(n x+3\beta)^2}{(1-\beta \log a)^3} + \frac{3 (n x+3\beta)\beta^2\log a }{(1-\beta \log a)^4} + \frac{ (\beta^3 + 2\beta^4)\log a  }{(1-\beta \log a)^5}\right)\right.\\
&&\left.+ 3\log a \left(\frac{n x+2\beta}{(1-\beta \log a)^2}+ \frac{\beta^2 \log a }{(1-\beta \log a)^3}\right) +\frac{1}{1- \beta \log a} \right\}\\
&=& \frac{x^3(\log a)^3 }{(1-\beta \log a)^3} + \frac{ 3 x^2 (\log a)^2 }{n (1-\beta \log a)^4} + \frac{x\log a(1+2 \beta  \log a +2 \beta ^4 (\log a)^3 -2 \beta ^4 (\log a)^4)}{n^2(1-\beta \log a)^5}.
\end{eqnarray*}
\begin{eqnarray*}
P_{n}^{[\beta,~a]}(t^4,x) &=& n x \sum_{k=0}^{\infty} (\log a)^k ( n x + k\beta)^{k-1} \frac{a^{-(nx+k\beta)}}{k!} \frac{k^4}{n^4}\\
&=&\frac{x}{n^3}\left\{ (\log a )^4S(4,n x+4\beta,\beta,a)+6(\log a )^3S(3,n x+3\beta,\beta,a)\right.\\
&&+7(\log a )^2S(2,n x+2\beta,\beta,a) +\log a ~ S(1,n x+\beta,\beta,a)\left.\right\}\\
&=&\frac{x}{n^3}\left\{\right. (\log a )^4\left(\frac{(n x+4\beta)^3 }{(1-\beta \log a)^4} +  \frac{6 (n x+4\beta)^2 \beta^2\log a }{(1-\beta \log a)^5 }+  \frac{2(n x+4\beta)\beta^3(2 +\beta)\log a+ 9(n x+4\beta) \beta^4(\log a)^2 }{(1-\beta \log a)^6 }\right.\\
&&\left.+ \frac{\beta^4 \log a +2\beta^5 (4  +\beta) (\log a)^2 + 4\beta^6(\log a)^3}{(1-\beta \log a)^7}\right)\\
&&+6(\log a )^3\left(\frac{(n x+ 3 \beta)^2}{(1-\beta \log a)^3} + \frac{3 (n x+ 3 \beta)\beta^2\log a }{(1-\beta \log a)^4} + \frac{ (\beta^3 + 2\beta^4)\log a  }{(1-\beta \log a)^5}\right)\\
&&+7(\log a )^2\left(\frac{(n x + 2 \beta)}{(1-\beta \log a)^2}+ \frac{\beta^2 \log a }{(1-\beta \log a)^3}\right)+\log a ~ \left(\frac{1}{1- \beta \log a}\right)\left.\right\}\\
&=&\frac{x^4 (\log a) ^4}{(1-\beta  \log a )^4} + \frac{6 x^3 (\log a) ^3}{n (1-\beta  \log a )^5} + \frac{x^2(\log a)^2 \left(7 +8 \beta  \log a+2 \beta ^4 (\log a) ^3-2 \beta ^4 (\log a) ^4\right)}{n^2 (1-\beta  (\log a) )^6}\\
&&+\frac{x \left((\log a) +8 \beta  (\log a) ^2+6 \beta ^2 (\log a) ^3+(12 \beta ^4 (\log a)^4-16 \beta ^5 (\log a) ^5+6 \beta ^6 (\log a) ^6)(1-\log a) \right)}{n^3 (1-\beta  (\log a) )^7}.
\end{eqnarray*}
This ends the proof.\\
We also introduce the $s$-th order central moment of the operator $P_{n}^{[\beta,~a]}$ , that is $P_{n}^{[\beta,~a]}(\varphi_x^s,x)$, where $\varphi_x(t) = t - x$, $(x, t) \in \mathbf{R}^+ \times\mathbf{R}^+$. On the basis of above lemma and by linearity of operators (\ref{14.1.eq14}), by a straightforward calculation, we obtain
\begin{lemma}\label{14.2.lemma4}
Let the operator $P_{n}^{[\beta,~a]}$  be defined by relation as (\ref{14.1.eq14}) and let $\varphi_x=t-x$ be given by
\begin{enumerate}
\item $\displaystyle P_{n}^{[\beta,~a]}(\varphi_x,x) = \frac{x (\log a +\beta  \log a -1)}{1-\beta  \log a}$;
\item $\displaystyle  P_{n}^{[\beta,~a]}(\varphi_x^2,x) = \frac{x \log a }{n (1-\beta  \log a )^3}+\frac{x^2 \left(1-2 (1+\beta ) \log a +(1+\beta )^2 (\log a)^2\right)}{(1-\beta  \log a )^2};$
\item $\displaystyle P_{n}^{[\beta,~a]}(\varphi_x^3,x) =\frac{x \left(\log a + 2 \beta  (\log a) ^2+ 2 \beta ^4 (\log a) ^4-2 \beta ^4 (\log a) ^5\right)}{n^2 (1-\beta  \log a )^5}+\frac{3 x^2\log a \left((1+\beta)  \log a-1\right)}{n (1-\beta  \log a )^4}\\
~~~~~~~~~~~~~~~~~~~~~~~~~~~~~~~~+\frac{x^3 \left(3 (\log a)  -(3 +6 \beta)  (\log a) ^2+(1+3 \beta+3 \beta ^2 )(\log a) ^3-(1-\beta \log a )^3\right)}{(1-\beta  \log a )^3};$
\item $\displaystyle P_{n}^{[\beta,~a]}(\varphi_x^4,x) =\frac{x \log a\left(1 + 8 \beta  \log a+ 6 \beta ^2 (\log a) ^2+ (12 \beta ^4 (\log a) ^3-16 \beta ^5 (\log a) ^4+6 \beta ^6 (\log a) ^5)(1- \log a)\right)}{n^3 (1-\beta  \log a )^7}\\
+ \frac{x^2 \log a\left((7-4 \beta ) \log a-4+8 \beta  (1+\beta ) (\log a) ^2-8 \beta ^4 (\log a) ^3+2 \beta ^4 (5+4 \beta ) (\log a) ^4-2 \beta ^4 (1+4 \beta ) (\log a) ^5\right)}{n^2 (1-\beta  \log a )^6}\\
+\frac{6 x^3\log a \left((1 -\beta\log a)^2-2 \log a+(1+2 \beta)  (\log a) ^2\right)}{n (1-\beta  \log a )^5}\\
+\frac{x^4\left(1-4 (1+\beta ) \log a +6 (1+\beta )^2 (\log a) ^2-4 (1+\beta )^3 (\log a) ^3+(1+\beta )^4 (\log a) ^4\right)}{(1-\beta  \log a )^4}. $
\end{enumerate}
\end{lemma}
\section{Approximation Properties}
The convergence property of the operator (\ref{14.1.eq14}) is proved in the following theorem:
\begin{theorem}
If $f\in  C[0, \infty)$ and $\beta_n \to 0,~a_n \to e $ as $n \to \infty$, then the sequence $P_{n}^{[\beta_n,~a]}$ converges uniformly to f(x) in $[a, b]$, where $0 \leq  a < b < \infty$.
\end{theorem}
\textbf{Proof: }Since $P_{n}^{[\beta_n,~a_n]}$ is a positive linear operator for $0\leq  \beta_n<1$ and $1< a_n \leq e$, it is
sufficient, by Korovkin's result \cite{korovkin1953convergence}, to verify the uniform convergence for test functions $f(t) = 1, t$ and $t^2$.\\
It is clear that $$P_{n}^{[\beta_n,~a_n]}(1,x) =1.$$
Going to $f(t) = t $,
$$\lim_{n\to \infty}P_{n}^{[\beta_n,~a_n]}(t,x) = \lim_{n\to \infty} \frac{x\log a_n}{1- \beta_n \log a_n}= x, \textrm { as } \beta_n \to 0 \textrm{ and  }a_n\to e.$$
Proceeding to the function $f(t) = t^2$, it can easily be shown that
$$\lim_{n\to \infty}P_{n}^{[\beta_n,~a_n]}(t^2,x) = \lim_{n\to \infty}\left(\frac{x^2(\log a_n)^2}{(1-\beta_n \log a_n)^2} +\frac{x \log a_n}{n(1-\beta_n  \log a_n)^3}\right) = x^2, \textrm { as } \beta_n \to 0 \textrm{ and  }a_n\to e,$$
and hence by Korovkin's theorem the proof of theorem is complete.\\
Following result gives order of approximation in term of modulus of continuity.
\begin{theorem}
If $f\in  C[0, \lambda]$ and $1 > \displaystyle \frac{\beta'}{n} \geq \beta \geq 0$ then
$$|f(x)- P_{n}^{[\beta,~a]}(f,x)| \leq \left\{ 1+  \left(\lambda\frac{1 +\lambda\beta\beta'}{(1-\beta  \log a )^3}\right)^{\frac{1}{2}}\right\}\omega\left(\frac{1}{\sqrt{n}}\right),$$
where $\omega(\delta) = \sup \{|f(x'')-f(x')|:~ x', x'' \in [0,\lambda]\}$, $\delta$ being a positive number such that $|x''-x'| < \delta$.
\end{theorem}
\textbf{Proof:}
By using the properties of modulus of continuity
\begin{equation}\label{14.2.eq3}|f(x'')-f(x')| \leq \omega(|x'-x''|); ~\textrm{ and }~ \omega(\gamma\delta) \leq (\gamma+1) \omega(\delta), ~~\gamma>0,\end{equation}
and noting the fact that
$$\sum_{k=0}^{\infty} \omega_{\beta,a}(k, n x)=1 \textrm{ and } \omega_{\beta,a}(k,nx )\geq 0 ~~~\forall~ n, k; 0\leq \beta<1 ; 1 < a\leq e.$$
It can easily be seen, by the application of Cauchy's inequality, that
\begin{eqnarray}\label{14.2.eq1}
|f(x)- P_{n}^{[\beta,~a]}(f,x)| &\leq& \left\{ 1+ \frac{1}{\delta}\sum_{k=0}^{\infty} n x (\log a)^k ( n x  + k\beta)^{k-1} \frac{a^{-(n x  +k\beta)}}{k!}\big |x-\frac{k}{n}\big| \right\}\omega(\delta)\nonumber\\
&\leq &\left\{ 1+ \frac{1}{\delta}\left[\sum_{k=0}^{\infty} n x (\log a)^k ( n x  + k\beta)^{k-1} \frac{a^{-(n x  +k\beta)}}{k!}\left(x-\frac{k}{n}\right)^2\right]^2 \right\}\omega(\delta).
\end{eqnarray}
Now by using second central moment (Lemma \ref{14.2.lemma4}), $0<\beta\leq 1$  and $1< a \leq e$, we have
\begin{eqnarray}\label{14.2.eq2}
\sum_{k=0}^{\infty} n x (\log a)^k ( n x  + k\beta)^{k-1} \frac{a^{-(n x  +k\beta)}}{k!}\left(x-\frac{k}{n}\right)^2 &=&  P_{n}^{[\beta,~a]}((t-x)^2,x)\leq \frac{x }{n (1-\beta  \log a )^3}+\frac{x^2\beta^2}{(1-\beta  \log a )^2}\nonumber\\
&\leq&  \frac{\lambda}{n}\left[\frac{1 }{(1-\beta  \log a )^3}+\frac{\lambda\beta\beta'}{(1-\beta  \log a )^2}\right]\nonumber\\
&\leq&  \frac{\lambda}{n}\left[\frac{1 +\lambda\beta\beta'}{(1-\beta  \log a )^3}\right].
\end{eqnarray}
Hence using (\ref{14.2.eq2}) in (\ref{14.2.eq1}) and choosing $\displaystyle \delta = \frac{1}{\sqrt{n}}$ we prove
$$|f(x)- P_{n}^{[\beta,~a]}(f,x)| \leq \left\{ 1+  \left(\lambda\frac{1 +\lambda\beta\beta'}{(1-\beta  \log a )^3}\right)^{\frac{1}{2}}\right\}\omega\left(\frac{1}{\sqrt{n}}\right).$$
For $\beta = 0$ and $a=e$, the above expression reduces to an inequality for the Szasz-Mirakyan operator obtained earlier by M\"{u}ller.
\begin{theorem}\label{14.2.thm2}
If $f\in C'[0, \lambda]$, $1 > \displaystyle \frac{\beta'}{n} \geq \beta \geq  0$, then the  following inequality holds
$$ |f(x)- P_{n}^{[\beta,~a]}(f,x)| \leq \left(\frac{\lambda(1 +\lambda\beta\beta')}{n(1-\beta  \log a )^3}\right)^{\frac{1}{2}}\left[1+\left(\frac{\lambda(1 +\lambda\beta\beta')}{n\delta(1-\beta  \log a )^3}\right)^{\frac{1}{2}}\right]\omega_1\left(\frac{1}{\sqrt{n}}\right). $$
where $\omega_1(\delta)$ is the modulus of continuity of $f'$.
\end{theorem}
\textbf{Proof:}
With out loss of generality we assume that $f'(x)\geq 0$. It may applies to $f'(x) < 0$. By the mean value theorem of differential calculus, it is known that
$$f(x) -f\left(\frac{k}{n}\right) = \left(x-\frac{k}{n}\right)f'({\xi}),$$
where $\xi = \xi_{n,k}(x)$ is an interior point of the interval determined by $x$ and $\displaystyle \frac{k}{n}$. Now
$$f(x) -f\left(\frac{k}{n}\right) \leq \left(x-\frac{k}{n}\right)[f'({\xi})-f'(x)] + \left[\frac{x\log a}{1-\beta \log a }-\frac{k}{n}\right]f'(x).$$
Multiplying both sides of the inequality by $\displaystyle nx(\log a)^k ( nx + k\beta)^{k-1} \frac{a^{-(nx+k\beta)}}{k!}$, summing
over $k$ and using $ P_{n}^{[\beta,~a]}(t,x)$, we have
\begin{equation}\label{14.2.eq4}|f(x)- P_{n}^{[\beta,~a]}(f,x)| \leq \sum_{k=0}^{\infty} \bigg|x-\frac{k}{n}\bigg|nx(\log a)^k ( nx + k\beta)^{k-1} \frac{a^{-(nx+k\beta)}}{k!}|f'({\xi})-f'(x)|.\end{equation}
But by (\ref{14.2.eq3}), we get
\begin{eqnarray*}
|f'({\xi})-f'(x)| \leq \omega_1(|\xi-x|) \leq \left(1+ \frac{1}{\delta}|\xi-x|\right)\omega_1(\delta)\leq \left(1+ \frac{1}{\delta}\bigg|\frac{k}{n}-x \bigg|\right)\omega_1(\delta),
\end{eqnarray*}
where $\delta$ is a positive number not depending on $k$.\\
A use of this in (\ref{14.2.eq4}) gives
\begin{eqnarray*}
|f(x)- P_{n}^{[\beta,~a]}(f,x)| &\leq& \left\{\sum_{k=0}^{\infty} \bigg|x-\frac{k}{n}\bigg|n x(\log a)^k ( n x + k\beta)^{k-1} \frac{a^{-(nx+k\beta)}}{k!}\right.\\
&&\left.+ \frac{1}{\delta}\sum_{k=0}^{\infty} \left(x-\frac{k}{n}\right)^2n x(\log a)^k ( n x + k\beta)^{k-1} \frac{a^{-(nx+k\beta)}}{k!}\right\}\omega_1(\delta).
\end{eqnarray*}
Hence by applications of Cauchy's inequality and (\ref{14.2.eq2}), we have
\begin{eqnarray*}
|f(x)- P_{n}^{[\beta,~a]}(f,x)| \leq \left(\frac{\lambda(1 +\lambda\beta\beta')}{n(1-\beta  \log a )^3}\right)^{\frac{1}{2}}\left[1+\left(\frac{\lambda(1 +\lambda\beta\beta')}{n\delta(1-\beta  \log a )^3}\right)^{\frac{1}{2}}\right]\omega_1(\delta).
\end{eqnarray*}
Choosing $\displaystyle \delta = \frac{l}{\sqrt{n}}$, Theorem \ref{14.2.thm2} is proved.\\
We know that a continuous function  $f$ defined on $I$ satisfies the condition  $$|f(x) -f(y)|\leq M_f|x-y|^{\eta}, (x,y) \in I \times E, $$ it called Locally $Lip~ \eta$ on $E(0< \eta \leq 1, E\subset I)$, where $M_f$ is a constant depending only on $f$.
\begin{theorem}
Let $E$ be any subset of $[0,\infty)$. If $f$ is locally $Lip~\eta$ on E, then we have
$$|P_{n}^{[\beta_n,~a_n]}(f,x) - f(x)| \leq M_fC(\eta,\beta,a)\max\{x^{\eta/2},x^{\eta}\} + 2 M_fd^{\eta}(x,E),$$
where $$\displaystyle C(\eta,\beta,a) = \left(\frac{1}{n (1-\beta  \log a )^3}+\frac{\beta^2}{(1-\beta  \log a )^2}\right)^{\eta}$$

and $d(x,E)$ is the distance between $x$ and $E$ defines as $$d(x,E) = \inf\{|x-y|:y\in E\}.$$
\end{theorem}
\textbf{Proof: }
Since $f$ is continuous,
$$|f(x) -f(y)|\leq M_f|x-y|^{\eta}$$ holds for any $x\geq 0$ and $y\in \overline{E}$, $\overline{E}$ is closure of $E \subset \mathbf{R}$. Let $(x,x_0) \in \mathbf{R}^{+} \times \overline{E}$  be such that $|x-x_0| = d(x,E)$.\\
Using Linear properties of $P_{n}^{[\beta_n,~a_n]}$, inequality  $(A+B)^{\eta}\leq A^{\eta} + B^{\eta}$ $( A\geq 0, B\geq 0, 0< \alpha \leq 1)$ and Holder's inequality , we get
\begin{eqnarray*}
P_{n}^{[\beta_n,~a_n]}(f,x) - f(x)| &\leq& P_{n}^{[\beta_n,~a_n]}(|f-f(x_0)|,x) + |f(x_0)- f(x)|\\
&\leq &  M_f P_{n}^{[\beta_n,~a_n]}(|(t-x)+(x-x_0)|^{\eta},x) + M_f |x_0- x|^{\eta}\\
&\leq& M_f \left(P_{n}^{[\beta_n,~a_n]}(|t-x|^{\eta},x) + |x-x_0|^{\eta}\right) + M_f|x_0- x|^{\eta}\\
&\leq& M_f \left(\left(P_{n}^{[\beta_n,~a_n]}((t-x)^{2},x)\right)^{\eta/2} +2 |x-x_0|^{\eta}\right) \\
&\leq& M_f \left(\left(\frac{x }{n (1-\beta  \log a )^3}+\frac{x^2\beta^2}{(1-\beta  \log a )^2}\right)^{\eta/2} + 2|x-x_0|^{\eta}\right) \\
&\leq & M_f \left(C(\eta,\beta,a)\max\{x^{\eta/2},x^{\eta}\} + 2 |x_0- x|^{\eta}\right)
\end{eqnarray*}
which is required results.\\
\section{Asymptotic Formula}
In order to present our asymptotic formula, we need the following lemma.
\begin{lemma}\label{14.3.lemma1} Let $ P_{n}^{[\beta_n,~a_n]}(f,x)$ be the generalized Jain operator. In addition, $\beta_n \to 0$ and $a_n \to e$ as $n \to \infty$, then
$$ P_{n}^{[\beta_n,~a_n]}((t-x)^4,x)\leq\frac{41(x + x^2+x^3+x^4) }{n^3(1-\beta_n \log a_n)^7}.$$
\end{lemma}
\textbf{Proof:}
Since $\max\{x, x^2,x^3,x^4\} \leq x + x^2+x^3+x^4$, $(1-\beta_n \log a_n)^2 \leq 1, \log a <1$ and $(1-\beta_n \log a_n)^{-i} \leq (1-\beta_n \log a_n)^{-(i+1)}$, for  $i=2,3,4,5,6$
\begin{eqnarray*}
P_{n}^{[\beta_n,~a_n]}((t-x)^4,x) &\leq& \frac{15 x}{n^3(1-\beta_n \log a_n)^7} + \frac{19 x^2}{n^2(1-\beta_n \log a_n)^6} + \frac{6 x^3}{n(1-\beta_n \log a_n)^5} + \frac{x^4}{(1-\beta_n \log a_n)^4}\\
&\leq& \frac{41(x + x^2+x^3+x^4) }{n^3(1-\beta_n \log a_n)^7}.
\end{eqnarray*}
\begin{theorem} Let $f,f',f'' \in  C(\mathbf{R}^+)$ and let the operator $P_{n}^{[\beta_n,~a_n]}$ be defined as in (1.9). If $\beta_n \to 0$ and $a_n \to e$ as $n \to \infty$ holds, then
$$\lim_{n\to \infty} n\left( P_{n}^{[\beta_n,~a_n]} (f,x) -f(x) \right) = \frac{x}{2} f''(x),~~ \forall ~~x>0. $$
\end{theorem}
\textbf{Proof: }
Let $f,f',f'' \in  C(\mathbf{R}^+)$ and $x \in [0,\infty)$ be fixed. By the
Taylor formula, we have
\begin{equation}\label{14.3.eq1}f(t) = f(x) + f'(x) (t-x) + \frac{1}{2}f''(x) (t-x)^2 + r(t;x)(t-x)^2,\end{equation}
where $r(t; x)$ is the Peano form of the remainder, $r(·; x) \in  C_2([0,\infty))$ and
$\displaystyle \lim_{t \to x} r(t,x) =0$.\\
We apply $P_{n}^{[\beta_n,~a_n]}$ to equation (\ref{14.3.eq1}), we get
\begin{eqnarray*}P_{n}^{[\beta_n,~a_n]}(f,x)- f(x)&=& f'(x) P_{n}^{[\beta_n,~a_n]}((t-x),x) + \frac{1}{2}f''(x) P_{n}^{[\beta_n,~a_n]}((t-x)^2,x) + P_{n}^{[\beta_n,~a_n]}(r(t ;x)(t-x)^2,x)\\
&=&  f'(x)\left[\frac{x\log a}{1- \beta \log a}-x \right]+ P_{n}^{[\beta_n,~a_n]}(r(t ;x)(t-x)^2,x) \\
&&+  \frac{f''(x)}{2}\left[\frac{x \log a }{n (1-\beta  \log a )^3}+\frac{x^2 \left(1-2 (1+\beta ) \log a +(1+\beta )^2 (\log a)^2\right)}{(1-\beta  \log a )^2}\right].\end{eqnarray*}
In the second  term $P_{n}^{[\beta_n,~a_n]}(r(t ;x)(t-x)^2,x)$ applying the Cauchy-Schwartz inequality, we have
\begin{equation}\label{14.3.eq2}0\leq |P_{n}^{[\beta_n,~a_n]}(r(t ;x)(t-x)^2,x)| \leq \sqrt{P_{n}^{[\beta_n,~a_n]}((t-x)^4,x)}\sqrt{P_{n}^{[\beta_n,~a_n]}(r(t ;x),x)}.\end{equation}
We have marked that $\displaystyle \lim_{t\to x}r(t, x) = 0$. In harmony with $\beta_n\to 0$ and $a_n\to e$ as $n\to \infty$, we have
\begin{equation}\label{14.3.eq3}\lim_{n\to \infty}  P_{n}^{[\beta_n,~a_n]}(r(t,x),x) =0.\end{equation}
On the basis of (\ref{14.3.eq2}), (\ref{14.3.eq3}) and Lemma \ref{14.3.lemma1} , we get that
$$\lim_{n\to \infty} n\left( P_{n}^{[\beta_n,~a_n]} (f,x) -f(x) \right) = \frac{x}{2} f''(x),~~ \forall ~~x>0. $$
\section{Statistical  Convergence}
Our next concern is the study of statistical convergence of the sequence of the generalized Sz\'{a}sz-Mirakyan operators .
we recall the concept of statistical convergence. This concept was introduced by Fast \cite{Fast1951sur}, The density of $S\subset \textbf{N}$ denoted by $\delta(S)$ is defined as follows
$$ \delta(S) = \lim_{n \to \infty} \frac{1}{n}\sum_{k=1}^n \chi_S(k),$$
provided the limit exists, where $\chi_S$ is the characteristic function of $S$.\\
A real sequence $x= (x_n)_{n\geq 1}$ is statistically convergent to real number $L$ if, for every $\epsilon>0$,
$$ \delta(\{ n\in \mathbf{N} : |x_n-L| \geq \epsilon \}) =0.$$
We write $\displaystyle st- \lim_n x_n =L$.  We note that every convergent sequence is a statistical convergent, but converse need not be true.
The application of this notion to the study of positive linear operators was first attempted in 2002 by Gadjiev and Orhan \cite{gadjiev2002some}.
\begin{theorem}\cite{gadjiev2002some}\label{14.4.thm1}
If $$st_A - \lim_n \|T_nx^k - x^k\| =0, k=0,1,2,$$
Then
$$st_A - \lim_n \|T_nf- f\| =0, $$
where $(T_n)_{n\geq 1}$ is a sequence of linear positive operators defined on $C(J)$ ; and the norm is taken on a compact interval included in $J$.
\end{theorem}
\begin{theorem}\label{14.4.thm2}
 Let $P_n^{[\beta_n,~a_n]}~,n\geq 1$ be defined as in (\ref{14.1.eq14}), where $(\beta_n)_{n> 1}, 0\leq \beta_n < 1$ and $1<a_n\leq e$, satisfies $\displaystyle st- \lim_n\beta_n =0 $ and $\displaystyle st- \lim_n a_n =e$.
Then
$$ st_A - \lim_n \|P_n^{[\beta_n,~a_n]} f - f \|_{[c,d]} = 0, ~~ f\in C(\mathbf{R^+}),$$
where $\|f\|_{[c,d]} = \sup_{x\in[c,d]} |f(x)|$.
\end{theorem}
\textbf{Proof: }
Notice that
$$\|P_n^{[\beta_n,~a_n]} (1,x) - 1 \|_{[c,d]} = 0, $$
$$\|P_n^{[\beta_n,~a_n]} (t,x) - x \|_{[c,d]} \leq \frac{d (\log a_n +\beta_n  \log a_n -1)}{1-\beta_n  \log a_n} , $$
$$\|P_n^{\beta_n,~a_n}( t^2,x) - x^2 \|_{[c,d]} \leq \frac{d \log a_n }{n (1-\beta_n  \log a_n )^3}+\frac{d^2 \left(1-2 (1+\beta_n ) \log a_n +(1+\beta_n )^2 (\log a_n)^2\right)}{(1-\beta_n  \log a_n )^2}. $$
Using $\displaystyle st_A- \lim_n \beta_n =0 $,  $\displaystyle st_A- \lim_n a_n =e $  and above inequality imply that,
$$st_A - \lim_n \|P_n^{[\beta_n,~a_n]}(t^k,x) -x^k\|_{[c,d]} = 0, ~~~k=0,1,2.$$
By theorem \ref{14.4.thm1}, we get required result.

\end{document}